\documentclass[preprint,12pt]{elsarticle}

\usepackage{graphicx}
\usepackage{amssymb}
 \usepackage{amsthm}

\usepackage{lineno}
\usepackage{amsmath}
\newtheorem{theorem}{Theorem}

\begin{document}
	\begin{center}
	\huge \bf	On the geometric structure of some statistical manifolds
	\end{center}
	
	\vspace*{1\baselineskip} 
	
	\begin{center}
		Mingao Yuan\\
			\vspace*{1\baselineskip} 
		\textit{Department of Statistics,}\\
		\textit{North Dakota State University,}\\
		\textit{Fargo, ND 58102, USA}\\
		E-mail: mingao.yuan@ndsu.edu
		
	\end{center}

\begin{center}	
\today	
\end{center}

	\vspace*{1\baselineskip} 
	
\begin{center}	
	\bf Abstract	
\end{center}

In information geometry, one of the basic problem is to study the geometric properties of statistical manifold. In this paper, we study the geometric structure of the generalized normal distribution manifold and show that it has constant $\alpha$-Gaussian curvature. Then for any positive integer $p$, we construct a $p$-dimensional statistical manifold that is $\alpha$-flat. 


	\section{Introduction}

	Fisher information is an important quantity in probability and statistics. It measures the amount of information that an observable random variable carries about the unknown parameters of the underlying distribution. The well-known Cramer-Rao theorem states that the lower bound of the variance of any unbiased estimator is the inverse of the Fisher information. In asymptotic theory, the maximum likelihood estimator converges in distribution to Gaussian distribution with mean zero and variance the inverse of the Fisher information. In 1945, Rao noticed that the Fisher information defines a Riemannian metric on a statistical manifold(\cite{rao}). Closely related to the Fisher information is the statistical curvature defined on one-parameter distribution family by Bradley Efron(\cite{efron}). It controls how much the variance of the maximum likelihood estimator exceeds the Cramer-Rao lower bound(\cite{efron}).  Later Madsen extended the result of Efron to the multi-parameter case(\cite{madsen}). It's well-known that differential geometry is an important field in mathematics. The famous Einstein's relativity theory depends on Riemannian geometry and recently some researchers are interested in extending the relativity theory by using the more general Riemann-Finsler geometry. See \cite{boothby,CS05, CY14,CZY14,SY16,YC15} for some references. 	
	In 1982, Amari provided a differential geometrical framework for analyzing statistical problmes related to mult-parameter families of distribution and introduced the $\alpha$-geometry on statistical manifold(\cite{amari1}). The $\alpha$-geometry measures the second-order information loss and second-order efficiency of an estimator(\cite{amari1}). Since then, many researchers studied the geometry of the statistical manifold(\cite{amari1}\cite{amari2}\cite{efron}\cite{kass}). Amari, Arwini and Dodson studied the $\alpha$-geometry of Gaussian, Gamma, Mckey bivariate gamma and the Freund bivariate exponential manifold(\cite{amari2}\cite{dodson}). Recently, the $\alpha$-geometry of Weibull, inverse gamma distribution, t-distribution and generalized exponential distribution manifold are investigated(\cite{caosun}\cite{cho}\cite{lisun}). One interesting fact is that the Gaussian manifold and the Weibull manifold have negative constant Gaussian curvature(\cite{dodson}\cite{caosun}) and several of the submanifolds of the Freund bivariate exponential manifold are $\alpha$-flat(\cite{dodson}). The statistical manifold with negative constant $\alpha$-curvature will share similar statistical properties as Gaussian manifold and Weibull manifold(\cite{amari1}\cite{efron}). Especially, the MLE for some parameter in $\alpha$-flat statistical manifold has no second order information loss(\cite{amari1}\cite{efron}). 
	 Then one both statistically and geometrically interesting question is whether we have other statistical manifolds that have constant Gaussian curvature or $\alpha$-Gaussian curvature. In this paper, we firstly show that the generalized Gaussian statistical manifold has constant $\alpha$-Gaussian curvature. Then for any positive integer $p$, we construct a $p$-dimensional statistical manifold that is $\alpha$-flat.

		The generalized Gaussian distribution is a generalization of the normal and Laplace distributions. It has received widespread applications in many applied areas(\cite{nada}\cite{pogany}). The generalized Gaussian distribution manifold is defined as
		\[M_1=\bigg \{f(x;\mu,\sigma,\beta)|f(x;\mu,\sigma,\beta)=\frac{\beta}{2\sigma\Gamma(1/\beta)}e^{-\frac{|x-\mu|^{\beta}}{\sigma^{\beta}}},\ x,\ \mu\in\mathbb{R},\ \sigma,\ \beta>0\bigg \},\]	
		where $\mu$, $\sigma$, $\beta$ are called the location, scale and shape parameters respectively and $\Gamma(x)$ is the gamma function. Clearly, this fimily includes the Gaussian distribution when $\beta=2$ and the Laplace distribution if $\beta=1$. Note that if $\beta$ is odd, the manifold is not smooth. Hence we only consider the case when $\beta$ is a known even number.

		\begin{theorem}
			Let $\beta$ be a given even number. Then the Riemannian metric on the generalized Gaussian statistical manifold $M_1$ is 
			\begin{equation}\label{gij}
			(g_{ij})=	 
			\begin{bmatrix}
			\frac{1}{\sigma^2}c_{11} & 0 \\
			0 & \frac{1}{\sigma^2}c_{22} 
			\end{bmatrix}	 ,
			\end{equation}	
			where 
			\[c_{11}=\frac{\Gamma(1-\frac{1}{\beta})\beta (\beta-1)}{\Gamma(\frac{1}{\beta})},\ \ \ \
				c_{22}=\beta.\]
			The $\alpha$-curvature tensor is given by
			\begin{equation}\label{rijkl}
		R^{(\alpha)}_{1212}=-\frac{(1-\alpha)\beta(\beta-1)[2-\beta+(1-\alpha)(\beta-1)]\Gamma(\frac{\beta-1}{\beta})}{\sigma^4\Gamma(\frac{1}{\beta})},
			\end{equation}
			and the $\alpha$-Gaussian curvature is constant and given by
			\begin{equation}\label{kcur}
			K^{(\alpha)}=-\frac{(1-\alpha)\big(2-\beta+(1-\alpha)(\beta-1)\big)}{\beta}.
			\end{equation}
			Then the $\alpha$-curvature tensor vanishes if and only if $\alpha=1$ or $\frac{1}{\beta-1}$.
		\end{theorem}

		Note that when $\alpha=0$, the $K^{(0)}$ is the Gaussian curvature of the Riemannian metric.	In this case, $K^{(0)}=-\frac{1}{\beta}$.
		If $\beta=2$, then the manifold is just the univariate Gaussian manifold. Plugging $\beta=2$ into the formula for $K^{(0)}$, we get $K^{(0)}=-\frac{1}{2}$, which is the same as that in (\cite{dodson}). When $\beta=4$, $K^{(0)}=-\frac{1}{4}$, and $K^{(0)}=-\frac{1}{6}$ for $\beta=6$. Then by using Theorem 1, we can get a lot of non-Gaussian statistical manifolds with constant Gaussian curvature different from that of the Gaussian manifold and the Weibull distribution.

	Next, we define another interesting non-Gaussian statistical manifold.	 
	 Let $\Omega_p=\{{\bf x}=(x_1, \dots, x_p)\in \mathbb{R}^p|\prod_{i=1}^px_i>0\}$ and $\mathbb{R}^p_+=\{{\bf x}=(x_1, \dots, x_p)\in \mathbb{R}^p|x_i>0, i=1, 2, \dots, p.\}$, we define a $p$-dimensional statistical manifold
	 
	 \[M_2=\bigg\{f({\bf x};{\bf \lambda})|f({\bf x};{\bf \lambda})=2\prod_{i=1}^p\frac{\sqrt{\lambda_i}}{\sqrt{2\pi}}e^{-\frac{\lambda_ix_i^2}{2}},\ {\bf x}\in\Omega_p,\ {\bf \lambda}\in\mathbb{R}^p_+ \bigg\}.\]
	 The importance of this distribution family lies in that its member is non-Gaussian multivariate distribution while the marginal distribution is Gaussian, which implies that a set of marginal distributions does not uniquely determine the multivariate normal distribution(\cite{dutta}).  For example, if $p=2$, we have 
	 	\begin{equation*}
	 	f(x_1, x_2)=2\frac{\sqrt{\lambda_1}}{\sqrt{2\pi}}e^{-\frac{\lambda_1x_1^2}{2}}\frac{\sqrt{\lambda_2}}{\sqrt{2\pi}}e^{-\frac{\lambda_2x_2^2}{2}}I[x_1x_2>0],
	 	\end{equation*}
	 	and the marginal distribution
	 	\begin{eqnarray*}
	 	f_{X_1}(x_1)&=&\int_{-\infty}^{+\infty}2\frac{\sqrt{\lambda_1}}{\sqrt{2\pi}}e^{-\frac{\lambda_1x_1^2}{2}}\frac{\sqrt{\lambda_2}}{\sqrt{2\pi}}e^{-\frac{\lambda_2x_2^2}{2}}I[x_1x_2>0]dx_2\\
	 	&=&\int_{-\infty}^{0}2\frac{\sqrt{\lambda_1}}{\sqrt{2\pi}}e^{-\frac{\lambda_1x_1^2}{2}}\frac{\sqrt{\lambda_2}}{\sqrt{2\pi}}e^{-\frac{\lambda_2x_2^2}{2}}I[x_1x_2>0]dx_2\\
	 	&&+\int_{0}^{+\infty}2\frac{\sqrt{\lambda_1}}{\sqrt{2\pi}}e^{-\frac{\lambda_1x_1^2}{2}}\frac{\sqrt{\lambda_2}}{\sqrt{2\pi}}e^{-\frac{\lambda_2x_2^2}{2}}I[x_1x_2>0]dx_2\\
	 	&=&\frac{\sqrt{\lambda_1}}{\sqrt{2\pi}}e^{-\frac{\lambda_1x_1^2}{2}}I[x_1<0]+\frac{\sqrt{\lambda_1}}{\sqrt{2\pi}}e^{-\frac{\lambda_1x_1^2}{2}}I[x_1>0]\\
	 	&=&\frac{\sqrt{\lambda_1}}{\sqrt{2\pi}}e^{-\frac{\lambda_1x_1^2}{2}},\ (x_1\in \mathbb{R}),
	 	\end{eqnarray*}
	 where $I$ is the indicator function.
	 Obviously $f$ is not a Gaussian density but $f_{X_1}$ is the density of the Gaussian distribution with mean zero and variance $\frac{1}{\lambda_1}$. Similarly, one can show that the another marginal distribution is Gaussian distribution with men zero and variance $\frac{1}{\lambda_2}$.

	  For this non-Gaussian manifold $M_2$, we have
	 
	 \begin{theorem}
	 	For any positive integer $p$, the $p$-dimensional statistical manifold $M_2$ is $\alpha$-flat.
	 \end{theorem}
	 
	 By this theorem, there exists $\alpha$-flat statistical manifold with any dimension.

	 \section{Geometry of statistical manifold}
	 
	 Let $M=\{p(x;\theta)|\theta\in\Theta\subset \mathbb{R}^p\}$ be a statistical manifold, $l=\log p(x;\theta)$ and $\partial_i=\frac{\partial}{\partial \theta_i}$. The Riemannian metric on $M$ is defined by
	\[g_{ij}(\theta)]=-E\big [\partial_i\partial_j l\big].\]
	The Levi-Civita connection is 
	 
	 		 \begin{equation*}
	 		 \Gamma^k_{ij}=g^{kl}\bigg\{\frac{\partial g_{li}}{\partial \theta_j}+\frac{\partial g_{lj}}{\partial \theta_i}-\frac{\partial g_{ij}}{\partial \theta_l}\bigg\},
	 		 \end{equation*}
	 		 and
	 \[\Gamma_{ijk}=\Gamma_{ij}^mg_{mk}.\]
	 
	 The $\alpha$-connection is defined by
	 \[\Gamma^{(\alpha)}_{ijk}=E\bigg[\bigg (\partial_i\partial_jl+\frac{1-\alpha}{2}\partial_il\partial_jl\bigg )\partial_kl\bigg].\]
	 Let $T_{ijk}=E\big[\partial_il\partial_jl\partial_kl\big]$ and $\Gamma^{(1)}_{ijk}=E[ \partial_i\partial_jl\partial_kl]$. Then we have
	 \[\Gamma^{(\alpha)}_{ijk}=\Gamma^{(1)}_{ijk}+\frac{1-\alpha}{2}T_{ijk}.\]

	Let $\Gamma^{(\alpha)k}_{ij}=g^{km}\Gamma^{(\alpha)}_{ijm}$. The $\alpha$-curvature tensor is 	
	 	\[R^{(\alpha)l}_{ihj}=\partial_i\Gamma^{(\alpha)l}_{hj}-\partial_h\Gamma^{(\alpha)l}_{ij}+\sum_m\Gamma^{(\alpha)l}_{im}\Gamma^{(\alpha)m}_{hj}-\sum_m\Gamma^{(\alpha)l}_{hm}\Gamma^{(\alpha)m}_{\ ij},\]
	 	and
	 	\[R^{(\alpha)}_{ihjk}=\sum_lg_{lk}R^{(\alpha)l}_{ihj}.\]
	 	
	 A statistical manifold is said to be $\alpha$-\textit{flat} if its $\alpha$-curvature vanishes.
	 For $p=2$, the $\alpha$-Gaussian curvature is defined as
	 \[K^{(\alpha)}=\frac{R^{(\alpha)}_{1212}}{det(g_{ij})}.\]

	Note that the $0$-geometry corresponds to the geometry of the Riemannian metric.

	 \section{Proof of the theorems}

	For distribution in $M_1$ and $\beta\neq1, 2$, we need to make a transformation of the parameter space so that the distribution can written as a regular exponential distribution. It is not easy to find such transformation for every $\beta$. So we work with the original parameter space without transformation. The computation is plausible.\\
	
	{\bf Proof of Theorem 1:} The log-likelihood function of the generalized normal distribution is
	\[l=\log f(x;\mu,\sigma,\beta)=\log\beta-\log\bigg (2\Gamma(\frac{1}{\beta})\bigg )-\log\sigma-\frac{(x-\mu)^{\beta}}{\sigma^{\beta}}.\]
	Then direct computation yields the first and second partial derivatives below
	\begin{eqnarray*}
	\frac{\partial l}{\partial \mu}&=&\frac{\beta}{\sigma^{\beta}}(x-\mu)^{\beta-1},\\
	\frac{\partial l}{\partial \sigma}&=&-\frac{1}{\sigma}+\frac{\beta}{\sigma^{\beta+1}}(x-\mu)^{\beta},\\
	\frac{\partial^2 l}{\partial \mu^2}&=&-\frac{\beta (\beta-1)}{\sigma^{\beta}}(x-\mu)^{\beta-2},\\
	\frac{\partial^2 l}{\partial \mu\partial \sigma}&=&	-\frac{\beta^2}{\sigma^{\beta+1}}(x-\mu)^{\beta-1},\\
	\frac{\partial^2 l}{\partial \sigma^2}&=&\frac{1}{\sigma^2}-\frac{\beta (\beta+1)}{\sigma^{\beta+2}}(x-\mu)^{\beta}.
	\end{eqnarray*}
	In terms of gamma function, we have the $k$-th moment
	\[  E[(x-\mu)^k]= \left\{
	\begin{array}{ll}
	0, & k: odd,\ \beta: even; \\
	\frac{\Gamma (\frac{k+1}{\beta})}{\Gamma(\frac{1}{\beta})}\sigma^k, & k,\ \beta: even.\\
	\end{array} 
	\right. \]
	Notice that we assume $\beta$ is even, then $\beta-1$ is odd. Hence, the Riemannian metric is
	\begin{eqnarray*}
	g_{11}&=&-E\bigg [\frac{\partial^2 l}{\partial \mu^2}\bigg ]=\frac{\beta (\beta-1)}{\sigma^{\beta}}E\bigg [(x-\mu)^{\beta-2}\bigg ]=\frac{\Gamma(1-\frac{1}{\beta})\beta (\beta-1)}{\Gamma(\frac{1}{\beta})}\frac{1}{\sigma^2},\\
	g_{22}&=&-E\bigg [\frac{\partial^2 l}{\partial \sigma^2}\bigg ]=-\frac{1}{\sigma^2}+\frac{\beta (\beta+1)}{\sigma^{\beta+2}}E\bigg [(x-\mu)^{\beta}\bigg ]\\
	&=&-\frac{1}{\sigma^2}+\frac{\beta (\beta+1)}{\sigma^{\beta+2}}\frac{\Gamma (\frac{\beta+1}{\beta})}{\Gamma(\frac{1}{\beta})}\sigma^{\beta}=\frac{\beta}{\sigma^2},\\
	g_{12}&=&g_{21}=-E\bigg [\frac{\partial^2 l}{\partial \mu\partial \sigma}\bigg ]=\frac{\beta^2}{\sigma^{\beta+1}}E\bigg [(x-\mu)^{\beta-1}\bigg ]=0,
	\end{eqnarray*}
	which leads to equation (\ref{gij}).
	
	Next we compute the coefficients $T_{ijk}$ below
	\begin{eqnarray*}
		T_{112}&=&E\bigg [\frac{\partial l}{\partial \mu}\frac{\partial l}{\partial \mu}\frac{\partial l}{\partial \sigma}\bigg]=-\frac{\beta^2}{\sigma^{2\beta+1}}E\bigg [(x-\mu)^{2(\beta-1)}\bigg ]+\frac{\beta^3}{\sigma^{3\beta+1}}E\bigg [(x-\mu)^{3\beta-2}\bigg ] \\
		&=&\frac{1}{\sigma^3}\frac{\Gamma(\frac{3\beta-1}{\beta})\beta^3-\Gamma(\frac{2\beta-1}{\beta})\beta^2}{\Gamma(\frac{1}{\beta})},\\
		T_{222}&=&E\bigg[\frac{\partial l}{\partial \sigma}\frac{\partial l}{\partial \sigma}\frac{\partial l}{\partial \sigma}\bigg]=	E\bigg[\bigg (-\frac{1}{\sigma}+\frac{\beta}{\sigma^{\beta+1}}(x-\mu)^{\beta}\bigg )^3\bigg]\\
		&=&-\frac{1}{\sigma^3}+\frac{3\beta}{\sigma^{\beta+3}}E\bigg [\bigg (x-\mu\bigg)^{\beta}\bigg]-\frac{3\beta^2}{\sigma^{2\beta+3}}E\bigg [\bigg (x-\mu\bigg)^{2\beta}\bigg]+\frac{\beta^3}{\sigma^{3\beta+3}}E\bigg [\bigg (x-\mu\bigg)^{3\beta}\bigg]\\
		&=&\frac{1}{\sigma^3}\bigg (-1+\frac{3\beta\Gamma(\frac{\beta+1}{\beta})}{\Gamma(\frac{1}{\beta})}-\frac{3\beta^2\Gamma(\frac{2\beta+1}{\beta})}{\Gamma(\frac{1}{\beta})}+\frac{\beta^3\Gamma(\frac{3\beta+1}{\beta})}{\Gamma(\frac{1}{\beta})}\bigg )\\
		&=&\frac{2\beta^2}{\sigma^3},\\
		T_{121}&=&T_{211}=T_{112},\\
		T_{111}&=&T_{221}=T_{212}=T_{122}=0.
	\end{eqnarray*}

	The 1-connection coefficients are
	
	\begin{eqnarray*}
	\Gamma^{(1)}_{112}&=&E\bigg [\frac{\partial^2 l}{\partial \mu^2}\frac{\partial l}{\partial \sigma}\bigg]=\frac{\beta (\beta-1)}{\sigma^{\beta+1}}E\bigg[(x-\mu)^{\beta-2}\bigg]-\frac{\beta^2 (\beta-1)}{\sigma^{2\beta+1}}E\bigg[(x-\mu)^{2\beta-2}\bigg],\\
	&=&\frac{1}{\sigma^3}\frac{\Gamma(\frac{\beta-1}{\beta})\beta (\beta-1)-\Gamma(\frac{2\beta-1}{\beta})\beta^2(\beta-1)}{\Gamma(\frac{1}{\beta})},\\
	\Gamma^{(1)}_{121}&=&E\bigg [\frac{\partial^2 l}{\partial \mu\partial \sigma}\frac{\partial l}{\partial \mu}\bigg]=-\frac{\beta^3}{\sigma^{2\beta+1}}E\bigg[(x-\mu)^{2\beta-2}\bigg]=-\frac{1}{\sigma^3}\frac{\Gamma(\frac{2\beta-1}{\beta})\beta^3 }{\Gamma(\frac{1}{\beta})},\\
	\end{eqnarray*}	
	\begin{eqnarray*}
	\Gamma^{(1)}_{222}&=&E\bigg [\frac{\partial^2 l}{\partial \sigma^2}\frac{\partial l}{\partial \sigma}\bigg]\\
	&=&-\frac{1}{\sigma^3}+\frac{\beta}{\sigma^{\beta+3}}E\bigg[(x-\mu)^{\beta}\bigg]+\frac{\beta (\beta+1)}{\sigma^{\beta+3}}E\bigg[(x-\mu)^{\beta}\bigg]-\frac{\beta^2(\beta+1)}{\sigma^{2\beta+3}}E\bigg[(x-\mu)^{2\beta}\bigg]\\
	&=&\frac{1}{\sigma^3}\bigg (-1+\frac{\beta\Gamma(\frac{\beta+1}{\beta})}{\Gamma(\frac{1}{\beta})}+\frac{\beta (\beta+1)\Gamma(\frac{\beta+1}{\beta})}{\Gamma(\frac{1}{\beta})}-\frac{\beta^2(\beta+1)\Gamma(\frac{2\beta+1}{\beta})}{\Gamma(\frac{1}{\beta})}\bigg )\\
		&=&-\frac{\beta (\beta-1)}{\sigma^3},
	\end{eqnarray*}	
	\begin{eqnarray*}
	\Gamma^{(1)}_{211}&=&\Gamma^{(1)}_{121},\\
	\Gamma^{(1)}_{111}&=&\Gamma^{(1)}_{122}=\Gamma^{(1)}_{212}=\Gamma^{(1)}_{221}=0.\\
	\end{eqnarray*}
	
	The $\alpha$-connection is just a linear combination of the $1$-connection and $T$. Hence, the $\alpha$-connection coefficients are
	
	\begin{eqnarray*}
	\Gamma^{(\alpha)}_{112}&=&\Gamma^{(1)}_{112}+\frac{1-\alpha}{2}T_{112}=\bigg (c^{(1)}_{112}+\frac{1-\alpha}{2}c_{112}\bigg )\frac{1}{\sigma^3},\\
	\Gamma^{(\alpha)}_{121}&=&\Gamma^{(1)}_{121}+\frac{1-\alpha}{2}T_{121}	=\bigg (c^{(1)}_{121}+\frac{1-\alpha}{2}c_{121}\bigg )\frac{1}{\sigma^3},\\
	\Gamma^{(\alpha)}_{222}&=&\Gamma^{(1)}_{222}+\frac{1-\alpha}{2}T_{222}	=\bigg (c^{(1)}_{222}+\frac{1-\alpha}{2}c_{222}\bigg )\frac{1}{\sigma^3},\\
	\Gamma^{(\alpha)}_{211}&=&\Gamma^{(\alpha)}_{121},\\
	\Gamma^{(\alpha)}_{111}&=&\Gamma^{(\alpha)}_{122}=\Gamma^{(\alpha)}_{212}=\Gamma^{(\alpha)}_{221}=0.\\	
	\end{eqnarray*}

	To compute the $\alpha$-curvature, we need the $\alpha$-connection coefficients in a another form.

	\begin{eqnarray*}
	\Gamma^{(\alpha)2}_{11}&=&g^{22}\Gamma^{(\alpha)}_{112}=\frac{1}{\sigma}\frac{1}{c_{22}}\bigg (c^{(1)}_{112}+\frac{1-\alpha}{2}c_{112}\bigg ),\\
	\Gamma^{(\alpha)1}_{21}&=&g^{11}\Gamma^{(\alpha)}_{211}=\frac{1}{\sigma}\frac{1}{c_{11}}\bigg (c^{(1)}_{121}+\frac{1-\alpha}{2}c_{121}\bigg ),\\
	\Gamma^{(\alpha)1}_{12}&=&\Gamma^{(\alpha)1}_{21},\\
	\Gamma^{(\alpha)2}_{21}&=&\Gamma^{(\alpha)2}_{12}=0.\\
	\end{eqnarray*}
	
	By definition, the $\alpha$-curvature is
	\begin{eqnarray}\nonumber
	R^{(\alpha)}_{1212}&=&-\Big[\bigg (\frac{\partial }{\partial \sigma}\Gamma^{(\alpha)2}_{11}-\frac{\partial }{\partial \mu}\Gamma^{(\alpha)2}_{21}\bigg )g_{22}+\Gamma^{(\alpha)}_{222}\Gamma^{(\alpha)2}_{11}-\Gamma^{(\alpha)}_{112}\Gamma^{(\alpha)1}_{21} \Big]\\   \label{new1}
	&=&-\frac{C_1+C_2-C_3}{\sigma^4},
	\end{eqnarray}
where the constants dependent on $\alpha$ and $\beta$ are defined below
			\begin{eqnarray*}
				c_{11}&=&\frac{\Gamma(1-\frac{1}{\beta})\beta (\beta-1)}{\Gamma(\frac{1}{\beta})},\\
				c_{22}&=&\beta,\\
				C_1&=&-\bigg (c^{(1)}_{112}+\frac{1-\alpha}{2}c_{112}\bigg ),\\
				C_2&=&\frac{1}{c_{22}}\bigg (c^{(1)}_{222}+\frac{1-\alpha}{2}c_{222}\bigg )\bigg ( c^{(1)}_{112}+\frac{1-\alpha}{2}c_{112}\bigg ),\\
				C_3&=&\frac{1}{c_{11}}\bigg (c^{(1)}_{112}+\frac{1-\alpha}{2}c_{112} \bigg )\bigg (c^{(1)}_{121}+\frac{1-\alpha}{2}c_{121} \bigg),\\
				\end{eqnarray*}
				\begin{eqnarray*}
				c_{112}&=&\frac{\Gamma(\frac{3\beta-1}{\beta})\beta^3-\Gamma(\frac{2\beta-1}{\beta})\beta^2}{\Gamma(\frac{1}{\beta})},\\
				c_{121}&=&c_{112},\\
				c_{222}&=&2\beta^2,\\
				c^{(1)}_{112}&=&\frac{\Gamma(\frac{\beta-1}{\beta})\beta (\beta-1)-\Gamma(\frac{2\beta-1}{\beta})\beta^2(\beta-1)}{\Gamma(\frac{1}{\beta})},\\
				c^{(1)}_{121}&=&-\frac{\Gamma(\frac{2\beta-1}{\beta})\beta^3 }{\Gamma(\frac{1}{\beta})},\\
				c^{(1)}_{222}&=&-\beta (\beta-1).\\
			\end{eqnarray*}	
	Then we can easily get the $\alpha$-Gaussian curvature below
	\begin{equation}\label{new2}
	K^{(\alpha)}=\frac{R^{(\alpha)}_{1212}}{det(g_{ij})}=-\frac{C_1+C_2-C_3}{c_{11}c_{22}}.
	\end{equation}
	
Next, we simplify (\ref{new1}) and (\ref{new2}), as pointed out by Professor Esmaeil Peyghan. Note that 
	\[C_1+C_2-C_3=\Big(c_{112}^{(1)}+\frac{1-\alpha}{2}c_{112}\Big)\Big(-1+\frac{c_{222}^{(1)}}{c_{22}}+\frac{1-\alpha}{2c_{22}}c_{222}-\frac{c_{121}^{(1)}}{c_{11}}-\frac{1-\alpha}{2c_{11}}c_{121}\Big).\]
	The first product factor can be calculated as
	\begin{eqnarray*}
	c_{112}^{(1)}+\frac{1-\alpha}{2}c_{112}&=&\frac{\Gamma(\frac{\beta-1}{\beta})\beta(\beta-1)[2-\beta+(1-\alpha)(\beta-1)]}{\Gamma(\frac{1}{\beta})},
	\end{eqnarray*}
	where we used the fact that $\Gamma(\frac{3\beta-1}{\beta})=\frac{(2\beta-1)(\beta-1)}{\beta^2}\Gamma(\frac{\beta-1}{\beta})$ and $\Gamma(\frac{2\beta-1}{\beta})=\frac{\beta-1}{\beta}\Gamma(\frac{\beta-1}{\beta})$, since $\Gamma(1+x)=x\Gamma(x)$.
	For the second product factor, we have 
	\begin{eqnarray*}
	-1+\frac{c_{222}^{(1)}}{c_{22}}-\frac{c_{121}^{(1)}}{c_{11}}&=&-\beta+\frac{\beta^2\Gamma(\frac{2\beta-1}{\beta})}{(\beta-1)\Gamma(\frac{\beta-1}{\beta})}
	=-\beta+\frac{\beta^2\frac{\beta-1}{\beta}\Gamma(\frac{\beta-1}{\beta})}{(\beta-1)\Gamma(\frac{\beta-1}{\beta})}=0,
	\end{eqnarray*}

	\begin{eqnarray*}
\frac{1-\alpha}{2c_{22}}c_{222}-\frac{1-\alpha}{2c_{11}}c_{121}&=&\frac{1-\alpha}{2}\Big(2\beta-\frac{\beta^2\Gamma(\frac{3\beta-1}{\beta})-\beta\Gamma(\frac{2\beta-1}{\beta})}{(\beta-1)\Gamma(\frac{\beta-1}{\beta})}\Big)\\
&=&\frac{1-\alpha}{2}\Big(2\beta-\frac{(2\beta-1)(
\beta-1)\Gamma(\frac{\beta-1}{\beta})-(\beta-1)\Gamma(\frac{\beta-1}{\beta})}{(\beta-1)\Gamma(\frac{\beta-1}{\beta})}\Big)\\
&=&1-\alpha.
	\end{eqnarray*}
	Then, we conclude that 
	\[	R^{(\alpha)}_{1212}=-(1-\alpha)\frac{\Gamma(\frac{\beta-1}{\beta})\beta(\beta-1)[2-\beta+(1-\alpha)(\beta-1)]}{\sigma^4\Gamma(\frac{1}{\beta})},\]
which is (\ref{rijkl}).	In this case, the $\alpha$-Gaussian curvature is

\begin{eqnarray*}
	K^{(\alpha)}&=&-(1-\alpha)\frac{\Gamma(\frac{\beta-1}{\beta})\beta(\beta-1)[2-\beta+(1-\alpha)(\beta-1)]}{\Gamma(\frac{1}{\beta})}\frac{\Gamma(\frac{1}{\beta})}{\Gamma(\frac{\beta-1}{\beta})\beta^2(\beta-1)}\\
	&=&-\frac{(1-\alpha)(2-\beta+(1-\alpha)(\beta-1))}{\beta},
\end{eqnarray*}
	which is (\ref{kcur}).

	\qed

		 \vspace*{1\baselineskip} 
	
	For distribution in $M_2$, we can easily write it as a regular exponential distribution. Then we can work on the potential function to get the $\alpha$-curvature(\cite{dodson}).\\
			
		{\bf Proof of Theorem 2:} We rewrite the distribution in $M_2$ as 
		\begin{eqnarray*}
			f({\bf x};{\bf \lambda})&=&e^{\frac{1}{2}\sum_{i=1}^p \log(\lambda_i)-\frac{1}{2}\sum_{i=1}^p\lambda_ix_i^2+\log 2-\log\sqrt{2\pi}}\\
			&=&e^{\frac{1}{2}\sum_{i=1}^p \log(-\theta_i)+\sum_{i=1}^p\theta_ix_i^2+\frac{p}{2}\log2-\log\sqrt{2\pi}},
		\end{eqnarray*}
		where $\theta_i=-\frac{1}{2}\lambda_i$. This is one member of the exponential family with $(\theta_1,\dots,\theta_p)$ the natural coordinates and the potential function
		\[\psi(\theta)=-\frac{1}{2}\sum_{i=1}^p \log(-\theta_i).\]
		
		For exponential family, the Fisher information is just the second derivative of the potential function(\cite{dodson}):
		\[g_{ij}=\frac{\partial^2\psi}{\partial\theta_i\partial\theta_j}=-\frac{1}{2}\frac{1}{\theta_i}\frac{1}{\theta_j}\delta_{ij}, \]
		where $\delta_{ii}=1$ for $i=1,\dots,p$ and $\delta_{ij}=0$ for $i\neq j$.
		The third derivative of the potential function will give us the $\alpha$-connection 
		\[\Gamma^{(\alpha)}_{ijk}=\frac{1-\alpha}{2}\frac{\partial^3\psi}{\partial\theta_i\partial\theta_j\partial\theta_k}=-\frac{1-\alpha}{2}\frac{1}{\theta_i}\frac{1}{\theta_j}\frac{1}{\theta_k}\delta_{ijk}, \]
		where $\delta_{iii}=1$ for $i=1,\dots,p$ and $\delta_{ijk}=0$ for unequal $i,j,k$.
		
		Then
		\[\Gamma^{(\alpha)k}_{ij}=g^{kl}\Gamma^{(\alpha)}_{ijl}=-\frac{1-\alpha}{(\theta_i\theta_j\theta_k)^{\frac{1}{3}}}\delta_{ijk}.\]

		Note that $\Gamma^{(\alpha)k}_{ij}$ and $\Gamma^{(\alpha)}_{ijk}$ vanish when $i,j,k$ are unequal.
		Hence the $\alpha$-curvature also vanish, that is,
		\[R^{(\alpha)}_{hijk}=0,\]
		which completes the proof.

		\qed \\

\section*{Acknowledgement}

Sincere thanks to Professor Esmaeil Peyghan for the valuable comments that significantly improve this manuscript.

\end{document}